\documentclass[12 pt]{article}
\usepackage{amsmath,amsthm,amssymb}
\usepackage[all]{xy}
\theoremstyle{plain}

\newtheorem{thm}{Theorem}[section]
\newtheorem{pro}[thm]{Proposition}
\newtheorem{lem}[thm]{Lemma}
\newtheorem{cor}[thm]{Corollary}
\newtheorem{con}[thm]{Conjecture}

\theoremstyle{definition}

\newtheorem{dfn}[thm]{Definition}
\newtheorem{rem}[thm]{Remark}

\DeclareMathOperator{\rk}{rk}
\DeclareMathOperator{\newmod}{mod}
\DeclareMathOperator{\card}{card}
\DeclareMathOperator{\res}{res}
\DeclareMathOperator{\codim}{codim}

\DeclareMathOperator{\FT}{FT}
\DeclareMathOperator{\divi}{div}
\DeclareMathOperator{\mult}{mult}

\DeclareMathOperator{\gr}{gr}
\DeclareMathOperator{\centre}{centre}
\DeclareMathOperator{\pr}{pr}
\DeclareMathOperator{\hcf}{hcf}

\DeclareMathOperator{\Proj}{Proj}
\DeclareMathOperator{\Diff}{Diff}
\DeclareMathOperator{\Exc}{Exc}
\DeclareMathOperator{\Hom}{Hom}
\DeclareMathOperator{\Spec}{Spec}

\renewcommand{\o}{\mathcal{O}}
\newcommand{\q}{\mathbb{Q}}
\newcommand{\R}{\mathbb{R}}
\newcommand{\p}{\mathbb{P}}
\newcommand{\C}{\mathbb{C}}
\newcommand{\z}{\mathbb{Z}}

\begin{document}
\bibliographystyle{halpha}

\title{Hypergeometric Equations and Weighted Projective Spaces}

\author{Alessio Corti\\
Department of Mathematics, Imperial College London\\
Huxley Building, 180 Queen's Gate\\
London SW7 2AZ, UK
\and Vasily Golyshev\\
Number Theory Section, Steklov Mathematical Institute\\
Gubkina str.\ 8,\\
119991, Moscow, Russia}

\date{30\textsuperscript{th} June, 2006}

\maketitle

\abstract{We compute the Hodge numbers of the polarised (pure)
  variation of Hodge structure $\mathbb{V}=\gr_{n-1}^W R^{n-1}f_{!}\,
  \mathbb{Z}$ of the Landau-Ginzburg model
  $f\colon Y \to \mathbb{C}$  mirror-dual to a weighted projective space
  $w\p^n$ in terms of a variant of Reid's age function of the
  anticanonical cone over
  $w\p^n$. This implies, for instance,
  that $w\p^n$ has canonical singularities if and
  only if $h^{n-1,0}\,\mathbb{V}=1$. We state a conjectural formula for
  the Hodge numbers of general hypergeometric variations.

  We show that a general fibre of the Landau-Ginzburg model is
  birational to a Calabi-Yau variety
  if and only if a general anticanonical section of $w\p$ is
  Calabi-Yau. We analyse the 104 weighted 3-spaces with canonical
  singularities,
  and show that a general anticanonical section is not a K3 surface
  exactly in those 9 cases where a generic fibre of the Landau-Ginzburg
  model is an elliptic surface of Kodaira dimension 1.}

\tableofcontents

\section{Introduction and results.}
\label{sec:introduction}

\subsection{Some hypergeometric local systems.}
\label{sec:main-result}

Fix positive integers (called weights) $w_0,\dots, w_n$ and write $d=\sum w_i$.
We assume that $\hcf (w_0, \ldots, \widehat{w_i},\ldots, w_n)=1$ for
all $i$, that is, the set of weights is \emph{well-formed}. The main
object of interest in this paper is the regular hypergeometric
operator on $\C^\times$:
\begin{multline}
\label{eq:1}
 H= \prod_{i=0}^n w_i^{w_i}D \Bigl(D-\frac1{w_i}\Bigr)\cdots
\Bigl(D-\frac{w_i-1}{w_i}\Bigr)-\\
-td^dD\Bigl(D+\frac1{d}\Bigr)\cdots\Bigl(D+\frac{d-1}{d}\Bigr)
\end{multline}
where $D=t\,d/dt$. Note that the symbol of $H$ is $\prod
w_i^{w_i}-td^d$; therefore, $H$ is singular at
$t=\lambda=(\prod w_i^{w_i})/(d^d)$.
More precisely, we are interested in the operator
$H^\text{red}$ obtained removing from each of the two summands on the
right hand side of Equation~\eqref{eq:1} one copy of every common factor;
it is known that $H^\text{red}$ is
irreducible.

In Theorem~\ref{thm:4}, we construct a precise Picard-Fuchs interpretation of
the local system of solutions of the ordinary differential equation
$H^\text{red}\varphi =0$ and, in Theorem~\ref{thm:1}, we compute its
Hodge numbers.

Consider the variety
\begin{equation}
  \label{eq:2}
  Y =
  \begin{cases}
    \prod_{i=0}^n y_i^{w_i} &=t\\
    \sum_{i=0} y_i & = 1
  \end{cases}
\qquad
\subset \quad \C^{\times \; n+1} \times \C^\times
\end{equation}
(where $y_0,\ldots,y_n$ and $t$ are coordinates on $\C^{\times \;
  n+1}$ and $\C^\times$) and the pencil $f=\pr_2 \colon Y \to \C^\times$.
\begin{thm}
  \label{thm:4}
  Denote by $\mathbb{V}_\R$ the real local system of solutions of the
  ordinary differential equation $H^\text{red}\varphi =0$\footnote{It
    is well-known from \cite{MR974906} that the (complex) local system
    of solutions of $H^\text{red}\varphi =0$ has a natural real
    structure.}; then
\[
\mathbb{V}_\R = \gr^W_{n-1} R^{n-1}f_!\, \R_Y.
\]
In particular, this naturally gives $\mathbb{V}_\R$ the structure of an integer
local system $\mathbb{V}_\z$.
\end{thm}

\begin{rem}
  \label{rem:2}
  Katz constructs a pure $\ell$-adic realisation of the local system
  $\mathbb{V}_\R$ \cite[Chapter 8]{MR1081536}; let us denote it by
  $\mathbb{V}_{\q_\ell}$;
  our proof of Theorem~\ref{thm:4} also shows that
  $\mathbb{V}_{\q_\ell} = \gr^W_{n-1} R^{n-1}f_!\, \q_{\ell\,Y}$.
\end{rem}

By Theorem~\ref{thm:4}, $\mathbb{V}_\z$ supports a polarised variation of
(integer) Hodge structure. Motivated by mirror symmetry, in
Theorem~\ref{thm:1}, we calculate the Hodge
numbers of this PVHS in terms of the geometry of the
\emph{anticanonical affine cone}
\[
A=\Spec \oplus_{m\geq 0} H^0(w\p,-mK)
\]
over the weighted projective space $w\p =\p(w_0,\dots,w_n)$.
Concretely, set $N=\z(w_0/d,\dots,w_n/d)+\z^{n+1}$ and let $M=\Hom
(N, \z)$; with this notation, $A$ is the affine toric variety
$\Spec \C\,[M\cap I^\ast]$
attached to the \emph{positive orthant} $I=\R_+^{n+1}
\subset N_\R$; equivalently, $A=\C^{n+1}/\mu_d$ is
the quotient of $\C^{n+1}$ by $\mu_d$ acting diagonally with weights
$w_0,\dots,w_n$.

We fix the linear form
\begin{equation}
  \label{eq:4}
  l\colon N \to \z \quad
  \text{defined as} \quad
  l(\alpha_0, \dots, \alpha_n)=\sum \alpha_i.
\end{equation}
If $a$ is an integer, we denote by $\overline a$ the smallest
positive integer $\newmod d$, that is, $0\leq \overline{a}<d$ and
$\overline{a} \equiv a \bmod d$.
Every element
\[
\Bigl(\frac{a_0}{d},\dots,\frac{a_n}{d}\Bigr) \in N
\quad
\text{has a unique representative}
\quad
\Bigl(\frac{\overline{a_0}}{d},\dots,\frac{\overline{a_n}}{d}\Bigr)
\]
$\bmod \, \z^{n+1}$ in the \emph{unit box} $[0,1)^{n+1} \subset \R^{n+1}$;
in particular, this identifies the cyclic group $\z/d= N/\z^{n+1}$
generated by $(w_0/d,\dots, w_n/d)$ with $N\cap [0, 1)^{n+1}$. We
denote by $a\colon N/\z^{n+1} \to [0,n]\cap \z$ the
\emph{age function} of Reid:
\[
a\Bigl(\frac{a_0}{d},\dots,\frac{a_n}{d}\Bigr)=
l\Bigl(\frac{\overline{a_0}}{d},\dots,\frac{\overline{a_n}}{d} \Bigr)
=
\frac1{d}\sum_{i=0}^n
\overline{a_i}.
\]
We denote by $(N/\z^{n+1})^0=(\z/d)^0$ the set of
classes not contained in a coordinate hyperplane:
\begin{multline*}
(\z/d)^0=\bigl\{0\leq k <d \mid d\; \text{divides
    no}\;kw_i\bigr\}=\\
   =\Bigl\{\Bigl(\frac{\overline{kw_0}}{d},\dots ,\frac{\overline{kw_n}}{d}
\Bigr) \mid \; \text{no} \; \overline{kw_i} = 0 \Bigr\}=  (N/\z^{n+1})^0,
\end{multline*}
and, for $j=1,\ldots,n$, we denote by
\[
a^s_j =  \card \, (\z/d)^0 \cap  a^{-1}(j)
\]
the number of elements of \emph{strict} age $j$; alternatively,
$a^s_j$ is the number
of elements $\mathbf{v} \in N \cap (0,1)^{n+1}$ with $l(\mathbf{v})=j$.
Our main result is:
\begin{thm}
  \label{thm:1} For $j=1, \ldots ,n$:
  \[
\rk \mathcal{H}^{n-j,j-1} \bigl(V_\z\otimes
\o_{\C^\times \smallsetminus \{\lambda \}}\bigr) =h^{n-j,j-1}\bigl(H^{n-1}_c
\,Y_t \bigr)= a^s_j.
  \]
\end{thm}
We prove Theorems~\ref{thm:4} and~\ref{thm:1} in \S\ref{sec:proof}. In
the remaining part of the Introduction, after a brief discussion of
some conjectures
and corollaries in the light of mirror symmetry, we
state our next results Theorems~\ref{thm:3}, \ref{thm:2} and~\ref{thm:5}.

\subsection{Hodge numbers of hypergeometric local systems}
\label{sec:hodge-numb-hyperg}

We expect that we will soon be able to announce progress on the
following natural conjecture:

\begin{con}
\label{conj_hodge_numbers}
For real numbers $0\leq \alpha_0 \leq\cdots \leq \alpha_{n-1}<1$ and
$0 \leq \beta_0 \leq\cdots \leq \beta_{n-1} < 1$ such that
$\alpha_j\not = \beta_k$ (all $j,k$) and the two sets $\{\exp 2\pi
i \alpha_j\}$, $\{\exp 2\pi i \beta_k \}$ are stable under complex
conjugation, consider the irreducible regular hypergeometric differential
operator:
\[ H\bigl(\{\alpha_j\};\{\beta_k\}\bigr)= \prod (D-\alpha_j)
-t \prod (D-\beta_k),
\]
write
\[
p(k)=\card \{ j \mid \alpha_j < \beta_k\} - k,
\]
and set
$p_+= \max \{p(k)\}, \;  p_-= \min \{p(k)\}$.
The local system of solutions of the ordinary differential
equation $H\varphi=0$ supports a real polarised variation of Hodge structure
of weight $p_+-p_-$ and Hodge numbers
\[h^{j-p_-,\,-j+p_+}=\card p^{-1} (j). \]
\end{con}

The main evidence for the conjecture is:

\begin{pro}
Given weights $w_0,\dots, w_n$ and $d=\sum w_i$ as before, set
\[
A = \bigcup_{i=0}^n
\Bigl\{\frac{k}{w_i} \mid k=0, \dots, w_i-1\Bigr\},
\quad
B = \Bigl\{\frac{k}{d} \mid k=0, \dots, d-1\Bigr\},
\]
and $\{ \alpha_j \} = A \setminus (A \cap B)$,
$\{ \beta_j \} = B \setminus (A \cap B)$.
Then the assertion of Conjecture \ref{conj_hodge_numbers} holds.
\end{pro}

\begin{proof}
Assume  $0 < k/d < 1$ is in
$\{ \beta _j \}$. This means precisely that
$k \in  (\z/d)^0$. One has
\[
p ( k ) = k - \sum_i \left( \left[ \frac{k w_i}{d} \right] + 1 \right),
\quad
\text{that is,}
\quad
p ( k ) = \left( \frac{\sum \overline{k w_i}}{d} \right) -
n.
\]
In other words, $p(k)=a(k)-n$; the
proposition now follows from Theorems \ref{thm:4} and \ref{thm:1}.
\end{proof}

Further evidence for the conjecture is given by the fact that it
correctly predicts the signature of the polarisation of a
hypergeometric local system computed by \cite[Theorem
4.5]{MR974906}. We are planning to cover the general case of the
conjecture in our next publication.

\subsection{Relation to the Reid-Tai criterion}
\label{sec:relation-reid-tai}

Let us look at the case $k=1$ of Theorem~\ref{thm:1}.
We view the affine hyperplane
$N_1 =N \cap l^{-1}(1)$ as a lattice by choosing the origin at
\[
\mathbf{e}=\Bigl(\frac{w_0}{d},\dots,\frac{w_n}{d}\Bigr).
\]
We regard $N_1=\Hom (\C^\times, \mathbb{T}^n)$ as the lattice of
one-parameter subgroups in a $n$-dimensional torus $\mathbb{T}^n$.
The primitive vectors
$\mathbf{e}_0,\dots ,\mathbf{e}_n \in N_1$
generate a simplex $\Delta$ containing the origin $\mathbf{e}$, and the
corresponding toric variety is the weighted projective space
$w\p=\p(w_0,\dots,w_n)$. The
Reid-Tai criterion \cite[\S4]{MR927963} states that the nonzero integer
points strictly inside $\Delta$ are in 1-to-1
correspondence with the geometric valuations $\nu$ of $w\p$ with
\emph{discrepancy} $a(\nu)<0$; therefore, we conclude:
\begin{cor}
  \label{cor:1}
  \[
h^{n-1,0}_c(Y_t)=1+\card \{\nu \mid a(\nu)<0\}.
  \]
In particular, $\p(w_0, \dots, w_n)$ has canonical singularities if and only if
$h^{n-1,0}_c(Y_t)=1$.
\end{cor}
(We recall the notion of discrepancy and canonical singularities in
Section~\ref{sec:canon-sing} below.)
\begin{rem}
  \label{rem:1}
  The linear function $l\colon N \to \z$ of Equation~\eqref{eq:4}
defines a grading on $\C[N]$ and a simplicial toric variety
$\p=\Proj \C[N\cap I]$ where $I=\R_+^{n+1}\subset N_\R$ is the positive
orthant; as we explain in the beginning of
Section~\ref{sec:proof-theorem-2}, the fibres $Y_t$ can naturally be seen to be
hypersurfaces in the torus $\mathbb{T}\subset \p$; therefore, they
are naturally compactified by simplicial hypersufaces
$\overline{Y_t}\subset \p$. The construction makes it clear that a
general $\overline{Y_t}$ is quasismooth; we conclude:
\[
h^{n-1,0}_c(Y_t)=h^0(\overline{Y_t},
K_{\overline{Y_t}})=h^0(\widetilde{Y_t}, K_{\widetilde{Y_t}})
\]
where $\widetilde{Y_t}$ is any nonsingular proper variety birational
to $Y_t$.
\end{rem}

\subsection{Relation to the McKay correspondence}
\label{sec:relat-mckay-corr}

The anticanonical affine cone $A=\Spec \oplus_{n\geq 0} H^0(w\p,-nK)$
over $w\p=\p(w_0, \dots, w_n)$ has, tautologically, Gorenstein
singularities; it has cyclic quotient singularities; therefore, it has
Gorenstein rational and hence canonical singularities.
The McKay correspondence \cite{MR1886756} states that, if
$\widetilde{A} \to A$ is a crepant resolution, then:
\[
b^{2j} (\widetilde{A})=\rk H^{2j}(\widetilde{A},\z)=
a_j =  \card \, (\z/d)\cap a^{-1}(j).
\]
More precisely, an element $g \in \z/d$ of age $j$ determines a
cohomology class in $H^{2j}(\widetilde{A},\z)$. We compare this
statement with our Theorem~\ref{thm:1}, stating that
\[
h^{n-j,j-1}\bigl(H^{n-1}_c
\,Y_t \bigr) = a^s_j =  \card \,(\z/d)^0 \cap a^{-1}(j).
\]
The comparison suggests:
\begin{con}
  $H\varphi=0$ is the quantum ordinary differential equation of the
  small quantum cohomology $QH^\bullet \widetilde{A}$, and
  $H^\text{red}\varphi =0$ is a
  direct summand of it. Moreover, $QH^\bullet \widetilde{A}$ is
  closely related to an appropriate variant of the small quantum orbifold
  cohomology $QH^\bullet_\text{orb} A$.
\end{con}
\emph{See} \cite{bryan-2005-} for the discussion of a very special
case. Note that $\widetilde{A}$ is non-compact and the cohomology ring
$H^\bullet\bigl(\widetilde{A}\bigr)$ does not satisfy Poincar\'e duality. Our
Theorem~\ref{thm:1} suggests that the subspace generated by the basis
elements corresponding to elements in $(\z/d)^0$ is a subring
satisfying Poincar\'e duality: does this subring have a natural
topological or geometric interpretation?

\subsection{Mirror symmetry}
\label{sec:mirror-symmetry}

We briefly discuss mirror symmetry for weighted projective spaces
$w\p=\p(w_0, \dots, w_n)$ in
the context relevant to this paper.
\cite{CCT} shows, in particular, that the small quantum orbifold cohomology
$QH^\bullet_\text{orb} (w\p)$ has rank $d=\sum w_i$ and, in the natural basis,
the quantum ordinary differential operator is
\[
P=-t+\prod_{i=0}^n w_iD\bigl(w_iD-1\bigr)\cdots \bigl(w_iD-(w_i-1)\bigr).
\]
If $H$ is the hypergeometric operator in Equation~\eqref{eq:1},
denote by $\mathcal{M}_{H^\text{red}}$ and $\mathcal{M}_P$ the
$\mathcal{D}$-modules on $\mathbb{C}^\times$ corresponding to
$H^\text{red}$ and $P$. The formulae in \cite[6.4.2]{MR1081536} state
that
\[
\FT j_{!\ast} [d]^\ast \mathcal{M}_P \cong j_{!\ast} [d]^\ast
\mathcal{M}_{H^\text{red}}
\quad \text{and} \quad
\FT j_{!\ast} [d]^\ast \mathcal{M}_{H^\text{red}} \cong j_{!\ast}
[d]^\ast
\mathcal{M}_P.
\]
In other words, the $d$-th Kummer pull-backs of the
$\mathcal{D}$-modules $\mathcal{M}_P$ and $\mathcal{M}_{H^\text{red}}$
are each other's Fourier transforms. This is a weak statement
of mirror symmetry between $w\p$ and the pencil $f\colon Y \to
\C^\times$. \emph{See} \cite{MR1874352}, where $\mathcal{M}_P$ and
$\mathcal{M}_{H^\text{red}}$ are called, respectively, the
\emph{Riemann-Roch} and the \emph{anticanonical Riemann-Roch}
$\mathcal{D}$-modules; the article \cite{math.AG/0404281} contains a much
deeper discussion of the case of weighted projective planes.

\subsection{Calabi-Yau weighted hypersurfaces}
\label{sec:calabi-yau-weighted}

The following notions are standard and well-known.

\begin{dfn}
  \label{dfn:2}
  \begin{enumerate}
  \item We say that a weighted projective space $w\p=\p(w_0, \dots, w_n)$ is
    \emph{well-formed} if $\hcf (w_0, \ldots, \widehat{w_i},\ldots,
     w_n)=1$ for all $i$.
  \item  A weighted hypersurface $X_e \subset w\p$ of degree $e$ is
    \emph{well-formed} \cite{MR1798982} if $X_e$ does not contain
    any codimension~2 coordinate subspace $(x_i=x_j=0)$.
  \item We say that $X_e$ is \emph{quasismooth} if it is well-formed
    and the affine cone over $X_e$ is nonsingular.
  \end{enumerate}
\end{dfn}

\begin{rem}
  \label{rem:3}
It is shown in \cite{MR1798982} that,
when $X_e$ is well-formed, the usual adjunction formula
$K_{X_e}=K_{w\p}+{X_e}_{|X_e}=\o_{X_e}(e-d)$ holds.
\end{rem}

\begin{dfn} \cite[Definition~4.1.8]{MR1269718}
  \label{dfn:3}
  A $n$-dimensional projective variety $Y$ is \emph{Calabi-Yau} if $Y$
  has canonical singularities, $K_Y\cong \o_Y$ (in particular, $Y$ is
  Gorenstein), and $H^i(Y, \o_Y)=(0)$ for $i=1, \dots, n-1$. A
  \emph{K3 surface} is a 2-dimensional Calabi-Yau variety.
\end{dfn}

(We recall the definition of canonical singularities in
Section~\ref{sec:canon-sing} below.)

\begin{lem}
  \label{lem:1}
  A quasismooth weighted hypersurface $X_d$ of degree $d=\sum w_i$
in $\p(w_0, \dots, w_n)$ is a Calabi-Yau variety.
\end{lem}

\begin{proof}
  Indeed, by Remark~\ref{rem:3}, $K_X \cong \o_X$; in particular, $X$
  has Gorenstein singularities. Because $X$ is quasismooth, it has
  cyclic quotient singularities; it follows that $X$ has Gorenstein
  rational, and hence canonical singularities. Finally, by
  \cite[Theorem~3.2.4]{MR704986}, $H^i(X, \o_X)=(0)$ for $i=1, \dots,
  n-1$.
\end{proof}

In the case of surfaces, a general
member $X_d\subset w\p^3$ is Calabi-Yau if and only if it is quasismooth:

\begin{thm} A general hypersurface
  \label{thm:3} $X_d \subset \p(w_0,\dots,w_3)$ is a well-formed
K3 surface if and only if it is quasismooth.
\end{thm}

Theorem~\ref{thm:3} is proved in Section~\ref{sec:proof-theorem-3}.

\begin{rem}
  \label{rem:4}
The above fails completely in higher dimensions.
For instance, in the list of Kreuzer and
Skarke\footnote{\emph{see} http:$/\!/$hep.itp.tuwien.ac.at/$\sim$kreuzer/CY/},
there are 184,026 Calabi-Yau weighted hypersurfaces in dimension 3;
of these, only 7,555 are quasismooth\footnote{\emph{see}
http:$/\!/$pcmat12.kent.ac.uk/grdb/}; a typical example is
$X_{99} \subset \p(15,18,19,20,27)$. (We are
grateful to M.~Reid for this remark.)
\end{rem}

The discussion above suggests that Calabi-Yau
hypersurfaces $X_d\subset w\p$ are mirror-dual to the fibres of the
pencil $f\colon Y \to \C^\times$ described by Equation~\eqref{eq:2}.
We investigate some crude geometric consequences of this statement and
we show that it is compatible with Batyrev's \cite{MR1269718} view of
mirror symmetry.

\begin{thm}
  \label{thm:2} Fix a well-formed weighted projective space
  $w\p=\p(w_0, \dots, w_n)$ and let $f\colon Y \to \C^\times$ be the
  pencil described by Equation~\eqref{eq:2}. A well-formed general hypersurface
  $X_d \subset w\p$ of degree $d=\sum w_i$ is
  Calabi-Yau if and only if a general fibre $Y_t=f^{-1}(t)$ is birational
  to a Calabi-Yau variety.
\end{thm}

\begin{rem}
  \label{rem:7}
  The proof of Theorem~\ref{thm:2}, which is carried out in
  \S\ref{sec:proof-theorem-3} below, is based on showing that,
  in the Calabi-Yau case, there are a resolution of $X_d$ and a
  nonsingular completion of $Y_t$ which are mirror-dual (nonsingular)
  Calabi-Yau toric
  hypersurfaces in the sense of \cite{MR1269718}.
\end{rem}

Presumably, it is possible to enumerate orbifold rational curves on a
quasismooth Calabi-Yau weighted hypersurface in terms of the polarised
variation of Hodge structure $\gr^W_{n-1}R^{n-1}f_!\, \z_Y$.

As we already remarked in Corollary~\ref{cor:1}, $w\p$ has canonical
singularities if and only if $h^{n-1,0}_c(Y_t)=1$. For $n=3$, there
are exactly $104$ weighted projective 3-spaces $\p(w_0,\dots,w_3)$
with canonical singularities (as can be checked by writing a small computer
program). The corresponding sets of weights are the \emph{famous 95}
of \cite[\S13.3, pp.\ 138--140]{MR1798982} plus the
\emph{additional nine} of Table~\ref{tab:weights}.
\begin{table}[ht]
  \centering
  \begin{tabular}[center]{cccc}
$w_0$ & $w_1$ & $w_2$ & $w_3$ \\
\hline
1 & 5 & 6 & 8 \\
\hline
1 & 4 & 7 & 9 \\
\hline
2 & 5 & 8 & 9 \\
\hline
1 & 5 & 8 & 14\\
\hline
3 & 7 & 8 & 10\\
\hline
4 & 7 & 9 & 10\\
\hline
5 & 8 & 9 & 11\\
\hline
3 & 7 & 8 & 18\\
\hline
5 & 8 & 9 & 22
  \end{tabular}
  \caption{The additional nine sets of weights}
  \label{tab:weights}
\end{table}
In the case of one of the additional nine sets of weights, a general surface
$X_d \subset w\p$ is well-formed but has non-canonical singularities;
nevertheless, $\gr^W_2 R^2f_! \z_Y$ is a polarised variation of Hodge
structure of weight 2 and $h^{2,0}=1$. It would be interesting to
study these PVHS in greater detail. In particular, it seems likely
that each fibre of the local system can be embedded in the K3 lattice and
hence, by Torelli, it is the transcendental lattice of a K3
surface. However, the local system as a whole may not be the PVHS
attached to the variation of transcendental lattice in a pencil of K3 surfaces,
possibly because the global monodromy acts nontrivially on the
discriminant group and hence does not lift to an action on the K3 lattice.

\begin{thm}
  \label{thm:5}
In the additional nine cases, $Y_t$ is birational to an elliptic surface
of Kodaira dimension $1$.
\end{thm}

Theorems~\ref{thm:3}, \ref{thm:2} and~\ref{thm:5} are proved in
Section~\ref{sec:proof-theorem-3}.

\section{Proof of Theorems~\ref{thm:4} and~\ref{thm:1}}
\label{sec:proof}

\paragraph*{Step 1}
\label{sec:step-1}
Recall the pencil of Equation~\eqref{eq:2}:
\[
  Y =
  \begin{cases}
    \prod_{i=0}^n y_i^{w_i} &=t\\
    \sum_{i=0} y_i & = 1
  \end{cases}
\qquad
\subset  \C^{\times \; n+1} \times \C^\times
\]
and $f=\text{pr}_2\colon Y \to \C^\times$. In
the notation of the Introduction:
\[
N=\z\Bigl(\frac{w_0}{d},\dots,\frac{w_n}{d}\Bigr)+\z^{n+1},
\]
$I =\R_+^{n+1}\subset N_\R$ is the positive
orthant, $l\colon N \to \z$ is the linear form of Equation~\eqref{eq:4},
$a\colon \z/d = N/\z^{n+1} \to [0,n]\cap \z$ is the age function,
\[
(\z/d)^0=(N/\z^{n+1})^0=\{0\leq k <d \mid d\; \text{divides no}\;kw_i\}
\]
and $a^s_j =  \card  a^{-1}(j)\cap (\z/d)^0$.

In this step, we show
that
\[
h^{n-j,j-1}\bigl(H^{n-1}_c
\,Y_t
\bigr)
= a^s_j.
\]

\begin{dfn}
  \label{dfn:8}
The \emph{Landau-Ginzburg pencil} or \emph{LG pencil} is the pull-back
$f^\prime \colon Y^\prime
\to \mathbb{C}^\times$
of $f\colon Y \to \C^\times$ by the covering $t=1/u^d\colon \C^\times
\to \C^\times$; after the change of coordinates $y_i\mapsto y_i/u$, we have:
\[
  Y^\prime = \Bigl(\prod_{i=0}^n y_i^{w_i} =1\Bigr) \subset \C^{\times \, n+1}
\quad
\text{and}
\quad
    f^\prime = u = \sum_{i=0}^n y_i.
\]
\end{dfn}
In what follows, it is convenient to work with the Landau-Ginzburg
pencil $u\colon Y^\prime \to \mathbb{C}^\times$; as in the
Introduction, we think of the affine hyperplane
$N_1 =l^{-1}(1)\cap N$ as a lattice by choosing
$(w_0/d,\dots,w_n/d)=\mathbf{e}$ as the origin; the projection
$\mathbf{w} \mapsto \mathbf{w}-l(\mathbf{w})\mathbf{e} +\mathbf{e}$
identifies $N_1$ with the lattice of characters of the $n$-dimensional torus
\[
\mathbb{T}^n=Y^\prime = \bigl\{(a_0,\dots, a_n)\in \C^{\times \, n+1}\mid \prod
a_i^{w_i}=1\bigr\}\subset \C^{n+1}.
\]
Denote by $\Delta=I\cap N_{1\,\R}$ the simplex generated by the vectors
$\mathbf{e}_i$. We may identify the Landau-Ginzburg fibre
$Y^\prime_u=Y_{1/t^d}$ with the hypersurface
\[
\bigl(-u\mathbf{e}+\sum \mathbf{e}_i =0\bigr)\subset \mathbb{T}^n=Y^\prime;
\]
in this description, it is clear that $\Delta$ is the Newton
polyhedron of $Y^\prime_u$, hence $Y^\prime_u$ is (tautologically)
Newton-regular (\emph{see} Definition~\ref{dfn:6}).

We use the notation and results of \cite[\S4 and \S5]{MR873655}; in
particular, we use the Poincar\'e series
\[
L_\Delta^\ast (t) = \sum_{k\geq 0} \ell^\ast (k\Delta)t^k
\]
where $\ell^\ast (k\Delta)$ is the number of lattice points strictly
inside $k\Delta$.
It is shown in \cite[\S4]{MR873655} that
\[
P_\Delta (t)=(1-t)^{n+1}L^\ast_\Delta (t)=\sum_{k=1}^n\varphi_k(\Delta)t^k
\]
is a polynomial. The statements \cite[3.11, 4.4, 5.6]{MR873655} imply that:
\begin{multline}
  \label{eq:3}
  h^{n-k, k-1}\bigl(H^{n-1}_c \,Y_t \bigr)=\\
= \varphi_k(\Delta)-\sum_{\dim \Gamma = n-1} \varphi_{k-1}(\Gamma)
+\sum_{\dim \Gamma = n-2} \varphi_{k-2}(\Gamma)-\cdots.
\end{multline}
We consider the locally closed
stratification of $I$ with strata corresponding to intersections of
coordinate hyperplanes and their translations by lattice vectors.
For $\mathbf{w} \in N\cap I$, we define the
\emph{content} $c(\mathbf{w})$:
\[
c(\mathbf{w})=k \quad \text{if $\mathbf{w}$ lies on the stratum of
  codimension $k$.}
\]
Let $B=\{b=(b_0,\dots,b_n)\subset \R^n \mid \text{all}\; 0\leq b_i <
1\}$ be the \emph{unit box}; note that:
\[
a_k^s=\card \{\mathbf{w} \in N\cap B \mid l(\mathbf{w})=k
\;\text{and}\;c(\mathbf{w})=0\}.
\]
We count using the exclusion-inclusion principle over strata of
codimension $k$ (the notation is self-explanatory):
\begin{multline*}
  \frac1{(1-t)^{n+1}}\sum_{k=1}^na_k^st^k=
  \sum_{\substack{\mathbf{w}\in N\cap
  I \\ c(\mathbf{w})=0}}t^{l(\mathbf{w})}=\\
=\sum_{\mathbf{w}\in N\cap I}t^{l(\mathbf{w})}-
\sum_{\substack{\mathbf{w}\in N\cap I \\ c(\mathbf{w})\geq 1}}
 t^{l(\mathbf{w})}=\\
=L^\ast_\Delta (t)-\frac{t}{1-t}\sum_{\codim F=1}L^\ast_F(t)+
\frac{t^2}{(1-t)^2}\sum_{\codim F =2} L^\ast_F(t)-\cdots
\end{multline*}
Therefore:
\begin{multline*}
  \sum a^s_k t^k = \\
=(1-t)^{n+1}L^\ast_\Delta(t)-t(1-t)^n\sum_{\codim F
  = 1} L^\ast_F(t)+\\
+t^2(1-t)^{n-1}\sum_{\codim F = 2}L^\ast_F
  (t)-\cdots =\\
= P_\Delta(t)-t\sum_{\codim F = 1}P_F(t)+t^2\sum_{\codim F = 2}
P_F(t)-\cdots=\\
=\sum_{k=1}^n h^{n-k,k-1}t^k.
\end{multline*}

\paragraph*{Step 2}
\label{sec:step-2}

The proof of \cite[Theorem~3.5.1]{MR1874352} shows that the
irreducible $\ell$-adic local system $\mathbb{V}_{\q_\ell}$ is contained
in $\gr^W_{n-1}R^{n-1}f_! \q_{\ell\, Y}$ for all $\ell$; it follows
that $\mathbb{V}_\R$ is a direct summand of $\gr^W_{n-1}R^{n-1}f_! \R_Y$;
summing the $h^{n-k, k-1}$ determined in Step~1 we calculate
the rank of $\gr^W_{n-1}R^{n-1}f_! \R_Y$ and conclude that the two
local systems have the same rank; therefore, they must
coincide. This proves Theorem~\ref{thm:4}, and the calculation in
Step~1 then proves Theorem~\ref{thm:1}. \qed

\section{Proof of Theorems~\ref{thm:3}, \ref{thm:2}, \ref{thm:5}}
\label{sec:proof-theorem-3}

\subsection{Toric varieties}
\label{sec:toric-varieties}

We use freely the language of toric varieties.
This section is here to fix our notation;
it is not an introduction to toric varieties.

\paragraph{Toric varieties}
\label{sec:toric-varieties-1}
If $\mathbb{T}$ is the $n$-dimensional torus, we denote by
$M=\Hom (\mathbb{T}, \mathbb{C}^\times)$ and
$N=\Hom (\mathbb{C}^\times, \mathbb{T})$ the
lattices of monomials (characters of $\mathbb{T}$) and weights
(one-parameter subgroups of $\mathbb{T}$).

We denote by $\p_\Sigma$ the toric variety associated to a
rational polyhedral fan $\Sigma\subset N_\R$. We mostly work with
complete fans, that is, proper toric varieties.

\paragraph{Weil and Cartier divisors and sheaves}
\label{sec:weil-cart-divis}
We denote by $\Sigma^{(1)}$ the set of primitive integer generators of the
1-dimensional cones of $\Sigma$; if $\mathbf{v}\in \Sigma^{(1)}$,
$D_\mathbf{v}\subset \p_\Sigma$ denotes the corresponding prime divisor;
for example, the canonical divisor is
\[
K=-\sum_{\mathbf{v}\in \Sigma^{(1)}} D_\mathbf{v}.
\]
Monomials $\mathbf{m}\in M$ are rational function on $\p_\Sigma$ and
\[
\divi \mathbf{m}=\sum_{\mathbf{v}\in \Sigma^{(1)}}
\langle \mathbf{m}, \mathbf{v}\rangle D_\mathbf{v}.
\]

An integer (rational) piecewise linear function $\varphi \colon N_\R \to \R$
which is linear on the cones of $\Sigma$ gives rise to a Cartier
($\q$-Cartier) divisor
\[
D=\sum_{\mathbf{v}\in \Sigma^{(1)}} \varphi(\mathbf{v}) D_\mathbf{v}
\]
on $\p_\Sigma$. The divisor $D$ is nef (ample) if and only if
$\varphi$ is convex (strictly convex). If $D= \sum_{\mathbf{v}\in \Sigma^{(1)}}
d_\mathbf{v}D_\mathbf{v}$ is a Weil divisor, we denote by
$\o(D)$ the sheaf of rational
functions $f$ such that $\divi f \geq -D$. When we view a rational
function $f\in \mathbb{C}[M]$ as a
rational section of $\o(D)$, the divisor of zeros and poles
of $f$ is $\divi f+D$. A monomial basis of
$H^0\bigl(\p_\Sigma, \o(D)\bigr)$ is
\[
\bigl\{\mathbf{m} \in M \mid \langle \mathbf{m},
\mathbf{v}\rangle \geq -d_\mathbf{v} \; \text{for all}\; \mathbf{v}\in
\Sigma^{(1)} \bigr\}\subset M.
\]

\paragraph{Polyhedra}
\label{sec:polyhedra}

The datum $(\Sigma, \varphi)$ of a fan $\Sigma$ and a
strictly convex integer piecewise linear function $\varphi\colon
N_\R \to \R$ linear on cones is equivalent to the datum of an integer
strictly convex polyhedron $Q\subset M_\R$ containing the
origin. Given $(\Sigma, \varphi)$, we set
\[
Q=\bigl\{\mathbf{m} \in M \mid \langle \mathbf{m},
\mathbf{v}\rangle \geq -\varphi(\mathbf{v}) \; \text{for all}\; \mathbf{v}\in
\Sigma^{(1)} \bigr\}\subset M.
\]

\paragraph{Toric hypersurfaces}
\label{sec:toric-hypersurfaces}

\begin{dfn}
  \label{dfn:6}
  A toric hypersurface $Z_f=(f=0)\subset X_\Sigma$ is
  \emph{Newton-regular} if it meets transversally all the toric strata.
\end{dfn}

\begin{rem}
  \label{rem:6}
  Consider a convex integer polyhedron $0 \in Q\subset M_\R$ and let
  $(\p_Q, L_Q)$ be the corresponding polarised toric variety.
  Consider a set of lattice points
  $S \subset Q\cap M$ and the linear system $\mathcal{D}=|S|
  \subset \p H^0(\p_Q, L_Q)$ that they generate. Then, a general
  member of $\mathcal{D}$ is Newton-regular if and only if $Q$ is the
  convex envelope of $S$.
\end{rem}

\subsection{Discrepancies and canonical singularities}
\label{sec:canon-sing}

We recall the definition of discrepancy of a valuation and the notion
of canonical singularities.

\begin{dfn}
  \label{dfn:1}
  \begin{enumerate}
  \item   Let $X$ be a normal variety of dimension $n$. We denote by $k(X)$
  the field of rational functions on $X$. A rational differential $\omega
  \in \Omega^n_{k(X)}$ determines a \emph{canonical divisor} of
  $X$:
\[
K_X =\divi_X \omega.
\]
  \item   A discrete rank 1 valuation $\nu$ of $k(X)$ is
  \emph{geometric} (sometimes called \emph{divisorial} in the
  valuation theory literature) if it has a \emph{uniformisation}, that
  is, a pair $E \subset Y$ of a normal variety $Y$ with $k(Y)=k(X)$,
  and a prime divisor $E\subset Y$ such that $\nu = \mult_E$ measures
  multiplicity along $E$. Abusing notation, we often make no
  distinction between the valuation $\nu$ and the divisor $E$.
  We say that $\nu$ has \emph{centre on $X$}, or that it is a
  valuation \emph{of} $X$ or \emph{on} $X$,
  if the obvious rational map $f\colon
  Y\dasharrow X$ is regular at $E$; the
  \emph{centre} of $\nu$ (or $E$) is the
  scheme-theoretic point $f(E)= \centre_X E$. We say that $\nu$ has
  \emph{small centre} if $\centre_X E$ has codimension $\geq 2$.
  \item Let $\nu$ be a geometric valuation with centre on
  $X$ and uniformisation $f \colon E\subset Y \to X$. If $K_X$ is
  $\q$-Cartier, we define the \emph{discrepancy} of $\nu$ as:
\[
a(\nu)=\mult_E \bigl(K_Y - f^\ast K_X\bigr)
\]
Note that, because $f$ is birational, $k(X)=k(Y)$; in the formula,
it is understood that $K_Y=\divi_Y \omega$. On the other hand, $f^\ast
K_X$ denotes the pull-back of a $\q$-Cartier divisor.
It is easy to see that the
discrepancy does not depend on the choice of uniformisation and the
differential $\omega$.
  \item Consider now a \emph{pair} $(X,B)$ of a normal variety $X$ and
  a \emph{boundary} $B$; that is, $B=\sum b_iB_i\subset X$ is a
  $\q$-divisor with prime components $B_i$ and $0<b_i\leq 1$. We say
  that $B$ is a \emph{strict} boundary if all $b_i<1$.
  If $\nu$ is a geometric valuation of $X$ with uniformisation
  $f\colon (E\subset Y)\to X$, and $K_X+B$ is $\q$-Cartier, we define
  the \emph{discrepancy} of $\nu$ with respect to the pair $(X,B)$ as:
\[
a(\nu, B)=\mult_E \bigl(K_Y - f^\ast (K_X+B)\bigr).
\]
  \end{enumerate}
\end{dfn}

\begin{dfn}
  \label{dfn:4}
  \begin{enumerate}
  \item $X$ has \emph{canonical} (\emph{terminal}) singularities if $K_X$ is
    $\q$-Cartier and every geometric valuation of $X$ has $\geq 0$ ($>0$)
    discrepancy.
  \item The pair $(X,B)$ has \emph{klt}, or \emph{Kawamata log
  terminal}, singularities if $K_X+B$ is
  $\q$-Cartier and every geometric valuation $\nu$ of $X$ has
  $a(\nu, B)>-1$.
  \item The pair $(X,B)$ has \emph{plt}, or \emph{purely log terminal},
  singularities if $K_X+B$ is
  $\q$-Cartier and every geometric valuation $\nu$ with small centre
  on $X$ has $a(\nu, B)>-1$.
  \item The pair $(X,B)$ has \emph{canonical}
  (\emph{terminal}) singularities if $K_X+B$ is
  $\q$-Cartier and every geometric valuation $\nu$ with small centre
  on $X$ has $a(\nu, B)\geq 0$ ($>0$).
  \item $X$ has \emph{pseudo-terminal} singularities if it has
  canonical singularities and $a(\nu)>0$ for every geometric valuation
  $\nu$ whose centre $\centre_X \nu$ has codimension $\geq 3$ in $X$.
  \end{enumerate}
\end{dfn}

\begin{rem}
  \label{rem:8}
  If $(X,B)$ has klt singularities, then $B$ is a strict boundary. It
  can be shown that if $(X,B)$ has plt singularities, then
\[
\lfloor B\rfloor = \sum_{b_i=1} B_i
\]
 is a disjoint union of normal components.
\end{rem}

\begin{dfn}
  \label{dfn:5}
  A \emph{good resolution} of a pair $(X, \mathcal{D})$ of a variety $X$
  and a linear system $\mathcal{D}$ on $X$ is a proper birational
  morphism $f\colon Y \to X$ where:
  \begin{enumerate}
  \item $Y$ is nonsingular and the exceptional set $\Exc f$ is a simple
  normal crossing divisor, and
  \item the birational transform $\mathcal{D}^\prime$ is a base point
  free linear system.
  \end{enumerate}
\end{dfn}

\begin{dfn}
  \label{dfn:7}
  If $X$ is a variety and $\mathcal{D}$ a linear system on $X$, we
  say that the pair $(X,\mathcal{D})$ has \emph{canonical}
  (\emph{terminal}) singularities if there is a good resolution
  $f\colon X \to Y$ with exceptional divisors $E_i\subset Y$ such that
\[
K_Y + \mathcal{D}^\prime =f^\ast (K_X+\mathcal{D})+\sum a_iE_i
\]
  with all $a_i\geq 0$ ($>0$).
\end{dfn}

\begin{rem}
  \label{rem:5}
  If the pair $(X,\mathcal{D})$ has \emph{canonical}
  (\emph{terminal}) singularities, and $D \in \mathcal{D}$ is a
  general member, then $(X,D)$ has canonical
  (terminal)  singularities. The converse is often not true.
\end{rem}

We refer the reader to \cite{MR927963} for an accessible introduction
to canonical singularities. It is easy to see that $X$ has
canonical singularities if $a(\nu)\geq 0$ for all $\nu$ which are
uniformised by a fixed resolution $f\colon Y \to X$.
When $\p$ is a toric variety, then $\p$ has a toric
resolution; in this case, it is enough to test toric valuations
associated to primitive vectors in the lattice $N$.

\subsection{Proof of Theorem~\ref{thm:2}}
\label{sec:proof-theorem}

Let $\mathbb{T}$ be the $n$-dimensional torus; denote by
$M=\Hom (\mathbb{T}, \mathbb{C}^\times)$ and
$N=\Hom (\mathbb{C}^\times, \mathbb{T})$ the
lattices of monomials and weights. The space $w\p^n=\p(w_0,\dots,
w_n)$ is a toric variety corresponding  to a choice of primitive
integer vectors $\mathbf{e}_i\in N$ satisfying the relation $\sum w_i
\mathbf{e}_i=\mathbf{0}$. In what follows, we denote by $\Delta
\subset N_\R$ the simplex generated by the $\mathbf{e}_i$.

Assume that a general $X_d\subset w\p^n$ is a well-formed Calabi-Yau.
Consider the dual simplex
\[\Delta^\ast=
\{\mathbf{m} \in M_\R \mid \text{all} \; \langle\mathbf{m},
\mathbf{e_i}\rangle \geq -1 \}\subset M_\R;
\]
note that $\Delta^\ast$ has rational not necessarily integer
vertices. The integer
points in $\Delta^\ast$ are a basis of $H^0(w\p, -K_{w\p})$. Consider
the integer polyhedron $Q\subset M$ generated by the integer points in
$\Delta^\ast$, and the associated polarised toric variety $(\p_Q, L_Q)$;
by construction, $L=-K_{\p_Q}$ and the proper transform of a general
$X_d \in |\o_{w\p}(d)|$ is a Newton-regular Calabi-Yau variety $X^\prime
\subset |-K_{\p_Q}|$.  By \cite[Theorem~4.1.9]{MR1269718}, $Q$ is a
reflexive polyhedron, and so is the dual polyhedron $P=Q^\ast \subset
N$. Denote by $\p_{P}$ the associated dual toric variety under the dual torus;
by construction, $\Delta \subset P$; therefore, the proper transform
\[
Y^{\prime\prime}_u=(-u\mathbf{0}+\sum_{i=0}^n \mathbf{e}_i= 0) \in |-K_{\p_P}|
\]
of the Landau-Ginzburg pencil is a pencil of anticanonical sections of
$\p_P$. The integer points in $P$ form a basis of
$H^0(\p_P, -K_{\p_P})$; more precisely,
an integer vector $\mathbf{v}\in P$ is a monomial on the dual torus; when
viewed as a section of $-K$, its divisor of zeros is
\[
Z(\mathbf{v})=\divi \mathbf{v} -K \sim -K
\]
(here $-K$ is the standard anticanonical divisor of $\p_P$,
that is, the sum of all the codimension 1 strata with multiplicity 1). Batyrev
\cite[Theorem~4.1.9]{MR1269718} states that a \emph{general} element of
$|-K_{\p_P}|$ is Calabi-Yau; we claim that a general $Y^{\prime\prime}_u$ is also
Calabi-Yau; in other words, we are saying that $Y^{\prime\prime}_u$ is
a projective birational Calabi-Yau model of the Landau-Ginzburg fibre
$Y^\prime_u$. This is the same as saying that a general member of the
linear system
\[
\mathcal{D} = |Z({\mathbf{e}_i}), Z(\mathbf{0})|
\]
is Calabi-Yau (absorb degrees of freedom using the torus
action); for this, it is enough to show that the pair
$(\p_P,\mathcal{D})$ has canonical singularities; that is, for all
valuations $\nu$ with small centre on
$\p_P$:
\[
a(\nu, \mathcal{D})=a(\nu, K_{\p_P})-\mult_\nu
\mathcal{D} \geq 0.
\]
It suffices to consider toric valuations corresponding to
primitive integer vectors $\mathbf{m} \in M\smallsetminus
\{\mathbf{0}\}$;
then:
\[
a(\mathbf{m},K_{\p_P})= -1-\mult_\mathbf{m} K
\]
and
\[
\mult_\mathbf{m} \mathcal{D} =-\mult_\mathbf{m} K
+ \min \{ \langle  \mathbf{m}, \mathbf{e}_i\rangle \};
\]
therefore, the statement follows from the claim that
\[
-1 \geq \min \{ \langle \mathbf{m},\mathbf{e}_i\rangle \}.
\]
This, however, is obvious: indeed, if it fails, that is precisely
saying that all $\langle \mathbf{m}, \mathbf{e}_i \rangle >-1$, that is,
$\mathbf{m}$ is strictly inside of $\Delta^\ast$. The assumption that
$X_d$ is a well-formed Calabi-Yau means exactly that $\mathbf{0}$ is
the only lattice point strictly inside of $\Delta^\ast$.

Vice-versa, assume that a general member $Y^\prime_u$ is birational to a
Calabi-Yau. We have already observed in the proof of
Theorem~\ref{thm:4} that $Y^\prime_u$ is Newton-regular.
By \cite{MR1690003}, we can
construct a (partial) resolution and minimal model $(\p_\Sigma,
Y^{\prime\prime}_u)$ where the birational transform
$Y^{\prime \prime}_u\in |-K_{\p_\Sigma}|$ is a Newton-regular
Calabi-Yau.  Again by \cite[Theorem~4.1.9]{MR1269718},
the toric variety $\p_\Sigma$ is the toric variety
$\p_P$ attached to a reflexive integer polyhedron $P\subset N$ and it is
obvious from the construction that $\Delta \subset P$. If
now $Q=P^\ast$ is the dual reflexive polyhedron, and $\p_Q$ is the
associated dual toric variety, the inclusion $\Delta \subset P$
gives a birational morphism $f\colon \p_Q \to w\p$.
By Batyrev, a general member $\widetilde{X} \in |-K_{\p_Q}|$ is a
Newton-regular Calabi-Yau; it follows that
$X=f(\widetilde{X})$ is Calabi-Yau.
\qed

\subsection{Proof of Theorem~\ref{thm:3}}
\label{sec:proof-theorem-1}

  In the proof, we need the following result of
  Shokurov and Koll\'ar \cite[Theorem~17.6]{MR1225842}.  Let $X$ a
  normal variety, $B=\sum b_i B_i$ a strict boundary, and $S\subset X$
  a normal subvariety of codimension 1 not contained in the support of
  $B$. As explained in \cite[Ch.\ 16]{MR1225842}, if $K+S+B$ is
  $\q$-Cartier, there is a naturally defined $\q$-divisor $B_S=\Diff_S
  B$ such that the adjunction formula
\[
(K_X+S+B)_S=K_S+B_S
\]
holds.

\begin{thm}[Inversion of adjunction]
  \label{thm:6}
  The pair $(X,S+B)$ has plt singularities in a neighbourhood of $S$ if
  and only if the pair $(S,B_S)$ has klt singularities. \qed
\end{thm}

\begin{proof}[Proof of Theorem~\ref{thm:3}] If
$X=X_d\subset w\p$ is a well-formed K3, then, by
Theorem~\ref{thm:6} with $S=X$ and $B=0$, the pair $(w\p, X)$ has plt
singularities hence ($K_{w\p}+X$ is Cartier, so all discrepancies
$a(E,X)$ are integers!) canonical
singularities in a neighbourhood of $X$; it follows that $w\p$ has
pseudo-terminal singularities along $X$; on the other hand, $w\p$ has
Gorenstein quotient hence canonical singularities outside of $X$;
therefore, $w\p$ has canonical singularities, hence it belongs to the
list of $104$.  We can check, by direct inspection, that if a general $X_d$ is
a K3, then it is a quasismooth K3: this is known in the
case of the 95; in the additional nine cases, we verify that $X_d$
is not a K3 surface.
\end{proof}

\subsection{Proof of Theorem~\ref{thm:5}}
\label{sec:proof-theorem-2}

The proof is an elementary verification. Recall the notation and
construction of the Landau-Ginzburg pencil:
\[
  Y^\prime = \Bigl(\prod_{i=0}^n y_i^{w_i} =1\Bigr) \subset \C^{\times \, n+1}
\quad
\text{and}
\quad
   u = \sum_{i=0} y_i\colon Y^\prime \to \C^\times.
\]
As usual we write $N=\z(w_0/d,\dots,w_n/d)+\z^{n+1}$ and $I
=\R_+^{n+1}\subset N_\R$; we denote by
$\mathbf{e}_i$ the standard basis vectors of $\z^{n+1}$ and
\[
\mathbf{e}=\Bigl(\frac{w_0}{d},\dots,\frac{w_n}{d}\Bigr)\in N.
\]
The linear function $l\colon N \to \z$ of Equation~\eqref{eq:4}
defines a grading on $\C[N]$ and the Landau-Ginzburg fibre
$Y^\prime_u\subset Y^\prime$
naturally compactifies to a hypersurface in a simplicial toric variety:
\[
\overline{Y^\prime_u}
=
\bigl(-u\mathbf{e}+\sum \mathbf{e}_i =0\bigr)\subset
\p=
\Proj \C[N\cap I].
\]
We use this compactification to study the Landau-Ginzburg pencil
$u\colon Y^\prime \to \C^\times$. To display information in
a human-readable format, it is convenient, though not necessary and
sometimes misleading, to realise the simplicial toric variety as a
quotient of a weighted projective space by the action of a product of
cyclic groups as in the following:

\begin{lem}
  \label{lem:2}
 Let $b_i=\hcf (w_i, d)$ and write $w_i=u_ib_i$, $d=d_ib_i$ with $\hcf
 (u_i, d_i)=1$. Then (in the notation of the preceding discussion)
\[
\p=\p(b_0, \dots, b_n)/G
\]
 where $G$ is a product of (at most) $n-1$ cyclic groups
acting diagonally and faithfully. Writing
$\tilde{\mathbf{e}}_i=\mathbf{e}_i/d_i$ and
$\widetilde{N}=\sum \z \tilde{\mathbf{e}}_i$, there are split exact
sequences:
\begin{gather*}
   0 \to \Hom \bigl(\widetilde{N}/N, \C^\times \bigr) \to  \prod \mu_{d_i}\to
 \mu_d \to 0,\\
 0 \to \mu_d \to \Hom \bigl(\widetilde{N}/N, \C^\times \bigr) \to G
 \to 0.
\end{gather*}
\end{lem}
\begin{proof}
  The proof is an easy exercise. It is clear that $\z^{n+1}\subset N
  \subset \widetilde{N}$ and that $\p$ is the quotient of
\[
\p(b_0, \dots, b_n)=\Proj \C[\widetilde{N}]
\]
 by the subgroup $\Hom \bigl(\widetilde{N}/N, \C^\times \bigr)$
 consisting of the
 elements of $\prod \mu_{d_i}= \Hom \bigl(\widetilde{N}/\z^{n+1},
 \C^\times \bigr)$ that fix
 $\mathbf{e}\in N$; this group does not act faithfully on $\p(b_0,
 \dots, b_n)$; dividing out by diagonal elements yields $G$.
\end{proof}

\begin{proof}[Proof of Theorem~\ref{thm:5}] The proof is a
  calculation summarised in Table~\ref{tab:ninewps}. Here, we explain
  what the table says and discuss the first entry in detail; the
  reader with sufficient motivation can perform
  the calculation in the remaining eight cases.

The first column of Table~\ref{tab:ninewps} lists the nine weights of
Table~\ref{tab:weights}; the second column gives the equation $F$ of
the pencil element $\overline{Y^\prime_u}=(F=0)$ in the natural
coordinates $z_0, \dots, z_3$ of $\p(b_0, \dots, b_3)$ corresponding
to the vectors $\tilde{\mathbf{e}}_0, \dots\tilde{\mathbf{e}}_3$
described in Lemma~\ref{lem:2}; the third column expresses the
simplicial toric variety $\p=\overline{Y^\prime}$ as $\p(b_0, \dots,
b_3)/G$ and identifies the group $G$; the fourth column gives the
unique $G$-invariant differential
\[
\omega = \res_{\,\overline{Y^\prime_u}} \, \frac{g(z) \Omega}{F}
\]
on $\overline{Y^\prime_u}$
where $\Omega=\sum (-1)^i z_i \, dz_0\cdots \widehat{dz_i}\cdots dz_3$
is the Griffiths differential on $\p(b_0,\dots, b_3)$; finally, the
last column lists a monomial $m(z_0, \dots, z_3)$ and the claim is
that $m\in N$ is of degree 0 and the rational map
\[
m\colon \overline{Y^\prime_u} \dasharrow \p^1
\]
is an elliptic fibration expressing $\overline{Y^\prime_u}$ as an
elliptic surface with Kodaira dimension $\kappa =1$.
\newcommand{\rb}[1]{\raisebox{1.5ex}[0pt]{#1}}
\begin{table}[ht]
\centering
  \begin{tabular}[center]{c|c|c|c|c}
    weights & $F$ & ambient & $\omega$ & map \\
    \hline
     & $z_0^{20}+z_1^4+z_2^{10}+z_3^5$ & & & \\
     \rb{1,5,6,8} &
     \hfill $-uz_0z_1z_2^{3}z_3^{2}$ &
     \rb{$\p(1,5,2,4)/\mu_{10}$} &
     \rb{$\frac{z_2^2z_3 \Omega}{F}$} &
     \rb{$\frac{z_3^2}{z_0^4z_2^2}$}\\
    \hline
     & $z_0^{21}+z_1^{21}+z_2^{3}+z_3^{7}$ & & & \\
     \rb{1,4,7,9} &
     \hfill $-uz_0z_1^{4}z_2z_3^{3}$ &
     \rb{$\p(1,1,7,3)/\mu_{21}$} &
     \rb{$\frac{z_1^{3}z_3^{2} \Omega}{F}$} &
     \rb{$\frac{z_{3}^3}{z_0^6z_1^3}$}\\
    \hline
    & $z_0^{12}+z_1^{24}+z_2^{3}+z_3^{8}$ & & & \\
     \rb{2,5,8,9} &
     \hfill $-uz_0z_1^{5}z_2z_3^{3}$ &
     \rb{$\p(2,1,8,3)/\mu_{12}$} &
     \rb{$\frac{z_1^{4}z_3^{2} \Omega}{F}$} &
     \rb{$\frac{z_{3}^3}{z_0^3z_1^3}$}\\
    \hline
    & $z_0^{28}+z_1^{28}+z_2^{7}+z_3^{2}$ & & & \\
     \rb{1,5,8,14} &
     \hfill $-uz_0z_1^{5}z_2^{2}z_3$ &
     \rb{$\p(1,1,4,14)/\mu_{14}$} &
     \rb{$\frac{z_1^{4}z_2\Omega}{F}$} &
     \rb{$\frac{z_{2}^2}{z_0^6z_1}$}\\
    \hline
    & $z_0^{28}+z_1^{4}+z_2^{7}+z_3^{14}$ & & & \\
     \rb{3,7,8,10} &
     \hfill $-uz_0^{3}z_1z_2^{2}z_3^{5}$ &
     \rb{$\p(1,7,4,2)/\mu_{14}$} &
     \rb{$\frac{z_0^{2}z_2z_3^{4} \Omega}{F}$} &
     \rb{$\frac{z_{2}^2}{z_0^4z_3^2}$}\\
    \hline
    & $z_0^{15}+z_1^{30}+z_2^{10}+z_3^{3}$ & & & \\
     \rb{4,7,9,10} &
     \hfill $-uz_0^{2}z_1^{7}z_2^{3}z_3$ &
     \rb{$\p(2,1,3,10)/\mu_{15}$} &
     \rb{$\frac{z_0z_1^6z_2^{2} \Omega}{F}$} &
     \rb{$\frac{z_2^2}{z_0^4z_1^4z_4}$}\\
    \hline
    & $z_0^{33}+z_1^{33}+z_2^{11}+z_3^{3}$ & & & \\
     \rb{5,8,9,11} &
     \hfill $-uz_0^{5}z_1^{8}z_2^{3}z_3$ &
     \rb{$\p(1,1,3,11)/\mu_{33}$} &
     \rb{$\frac{z_0^4z_1^7z_2^2 \Omega}{F}$} &
     \rb{$\frac{z_2^3}{z_0^6z_1^3}$}\\
    \hline
    & $z_0^{12}+z_1^{36}+z_2^{9}+z_3^{2}$ & & & \\
     \rb{3,7,8,18} &
     \hfill $-uz_0z_1^7z_2^{2}z_3$ &
     \rb{$\p(3,1,4,18)/\mu_{6}$} &
     \rb{$\frac{z_1^{6}z_{2} \Omega}{F}$} &
     \rb{$\frac{z_{2}^2}{z_0^2z_1^2}$}\\
    \hline
    & $z_0^{44}+z_1^{11}+z_2^{44}+z_3^{2}$ & & & \\
     \rb{5,8,9,22} &
     \hfill $-uz_0^{5}z_1^{2}z_2^{9}z_3$ &
     \rb{$\p(1,4,1,22)/\mu_{22}$} &
     \rb{$\frac{z_0^{4}z_1z_2^{8} \Omega}{F}$} &
     \rb{$\frac{z_{1}^2}{z_0^6z_2^2}$}\\
    \hline
  \end{tabular}
  \caption{Proof of Theorem~\ref{thm:5}}
  \label{tab:ninewps}
\end{table}
We now verify all these statements for the first set of weights
$(w_0,w_1,w_2, w_3)=(1,5,6,8)$. We leave the verification in the other eight
cases to the reader.

According to Lemma~\ref{lem:2}, we have
\[
\widetilde{N}=\tilde{\mathbf{e}}_0\z+\tilde{\mathbf{e}}_1\z+
\tilde{\mathbf{e}}_2\z+\tilde{\mathbf{e}}_3\z=
\frac{\mathbf{e_0}}{20}\z+\frac{\mathbf{e}_1}{4}\z+
\frac{\mathbf{e}_2}{10}\z+\frac{\mathbf{e}_3}{5}\z,
\]
and
\[
\mathbf{e}=\frac{1}{20}(1,5,6,8)=\frac{\mathbf{e}_0}{20}+
\frac{\mathbf{e}_1}{4}+3\frac{\mathbf{e}_2}{10}+
2\frac{\mathbf{e}_3}{5};
\]
this shows that
\[
\overline{Y^\prime_u}=
\bigl(z_0^{20}+z_1^4+z_2^{10}+z_3^5-uz_0z_1z_2^3z_3^2=0\bigr)\subset
\p(1,5,2,4)/G.
\]
It is easy to see that $G=\mu_{10}$ acts on $\p(1,5,4,2)$
with weights $0,5,1,6$:
\[
z_0, z_1, z_2, z_3 \mapsto z_0, \zeta^5 z_1, \zeta z_2, \zeta^6 z_3
\]
where $\zeta$ is a primitive $10$-th root of unity.

Note that
\[
\frac{\Omega}{z_0z_1z_2z_3}
\]
is a $\mathbb{T}$-invariant, hence $G$-invariant, rational
differential on $\p(1,5,2,4)$; it follows that
\[
\omega = \res_{\, \overline{Y^\prime_u}} \, \Bigl( z_0z_1z_2^3z_3^2
\frac{\Omega}{z_0z_1z_2z_3 F}\Bigr)
\]
is a $G$-invariant regular differential on $\overline{Y^\prime_u}$. This
differential must be unique since, after all, the list was produced by imposing
$h^{2,0}(\overline{Y^\prime_u})=h^{2,0}_c(Y^\prime_u)=1$,
cf.\ Remark~\ref{rem:1}.

Now we claim that the curve $C=(z_3=0)\subset \overline{Y^\prime_u}$ is
a nonsingular curve of genus 1. This is a completely elementary
calculation. Indeed, consider the lattice
\[
N_{z_3} = N\cap (\alpha_3=0)=\mathbf{e}_0\z+\mathbf{e}_1\z
+\mathbf{e}_2\z+\frac{1}{4}(1,1,2,0)\z;
\]
$C$ naturally sits in the simplicial plane $\Proj \C [N_{z_3}]$;
using Lemma~\ref{lem:2} again, we can write
\[
C=(x_0^4+x_1^4+x_2^2=0)\subset \p(1,1,2)/\mu_{2}
\]
where $\mu_2$ acts diagonally as $(x_0,x_1,x_2)\mapsto (x_0, -x_1,
-x_2)$. The map to the quotient is \'etale over $C$; this description
hence shows that $C$ is a nonsingular curve of genus 1. This proves,
in particular, that the Kodaira dimension of $\overline{Y^\prime_u}$
is $\geq 1$; indeed, the curve C is contained in the zero locus
$(z_2^2z_3=0)=Z(\omega)$ of the differential $\omega$ and, because $C$
is not a rational curve, $C$ survives to give a nontrivial element
in the canonical class of the minimal model of
$\overline{Y^\prime_u}$.

To conclude, it remains to show that the map
\[
m=\frac{z_3^2}{z_0^4 z_2^2}\colon \overline{Y^\prime_u} \dasharrow \p^1
\]
is an elliptic fibration. First of all, we observe that the map is
well defined, that is, the monomial defining it is $G$-invariant or,
which is the same, it is an element of $N$ of degree 0. But indeed
\[
\frac{z_3^2}{z_0^4 z_2^2}=\frac{z_0z_1z_2^3z_3^2}{z_0^5z_1z_2^5}
\quad
\text{is the vector}
\quad
\mathbf{e}-\frac{1}{4}(1,1,2,0) \in N
\]
of degree $1-1=0$. Thus, $m \colon \overline{Y^\prime_u}\dasharrow \p^1$ is a
well defined map.
The vectors $x_i=\mathbf{e}_i$, $y=\mathbf{e} \in \mathbb{C}[N\cap I]_1$
define a birational morphism $\pi \colon \p\to \p^4$ with image
the non-normal hypersurface
\[
y^{20}-x_0x_1^5x_2^6x_3^8=0.
\]
The image $\pi(\overline{Y^\prime_u})$ is the complete intersection:
\begin{equation}
  \label{eq:5}
  \begin{cases}
  y^{20}-x_0x_1^5x_2^6x_3^8=0\\
  -uy+x_0+x_1+x_2+x_3=0
\end{cases}
\subset \quad \p^4.
\end{equation}
Note that
\[
\frac{z_3^2}{z_0^4z_2^2}=\frac{z_0z_1z_2^3z_3^2}{z_0^5z_1z_5}=
\frac{x_1x_2x_3^2}{y^4};
\]
therefore, omitting coefficients, and substituting $y^4=x_1x_2x_3^2$ in
the first line of Equation~\eqref{eq:5}, we need to check that
\begin{equation}
  \label{eq:6}
  \begin{cases}
  x_3^2-x_0x_2=0\\
  -y+x_0+x_1+x_2+x_3=0\\
  y^4-x_1x_2x_3^2=0
\end{cases}
\subset \quad \p^4
\end{equation}
is birational to a curve of genus 1. This is an elementary
verification: eliminating $x_1$ and then $x_0$ we obtain
\[
y^2+y(x_3^2 +x_3x_4)+x_3^4+x_4^4=0
\]
which is, manifestly, the equation of a curve of genus 1.
\end{proof}

\addcontentsline{toc}{section}{References}

\bibliography{cg.bib}

\end{document}